\documentclass[12pt,reqno]{amsart} 
\usepackage{amsthm,amsfonts,amssymb,amscd}

\newtheorem{theorem}[subsection]{Theorem}
\newtheorem{proposition}[subsection]{Proposition}
\newtheorem{prop}[subsection]{Proposition}

\newtheorem{lemma}[subsection]{Lemma}
\newtheorem{corollary}[subsection]{Corollary}

\theoremstyle{definition}
\newtheorem{definition}[subsection]{Definition}
\newtheorem{example}[subsection]{Example}
\newtheorem{proposition-definition}[subsection]{Proposition-Definition}

\theoremstyle{remark}
\newtheorem{remark}[subsection]{Remark}

\newcommand{\mt}[1]{\operatorname{#1}}

\newcommand{\Diff}{\operatorname{Diff}}
\newcommand{\Supp}{\operatorname{Supp}}
\newcommand{\mult}{\operatorname{mult}}

\newcommand{\pal}{\text{---}}

\newcommand{\GL}{{GL}_3(\CC)}
\newcommand{\SL}{{SL}_3(\CC)}
\newcommand{\PGL}{{PGL}_3(\CC)}
\newcommand{\red}{\operatorname{red}}

\newcommand{\CC}{{\mathbb C}}
\newcommand{\RR}{{\mathbb R}}
\newcommand{\ZZ}{{\mathbb Z}}
\newcommand{\QQ}{{\mathbb Q}}
\newcommand{\PP}{{\mathbb P}}

\newcommand{\CCC}{{\mathcal C}}
\newcommand{\KKK}{{\mathcal K}}

\newcommand\alp{\alpha}
\newcommand\de{\delta}
\newcommand\De{\Delta}

\newcommand\la{\lambda}

\newcommand\G{\Gamma}

\newcommand\lra{{\longrightarrow}}

\author{D. Markushevich}

\address{D. M.: Math\'ematiques - b\^{a}t. M2, Universit\'e Lille 1,
F-59655 Villeneuve d'Ascq Cedex, France}
\email{markushe@gat.univ-lille1.fr}

\author{Yu.~G.~Prokhorov}

\thanks{The second author was partially supported
by the Russian Foundation of Fundamental Research grant
96-01-00820 and by a grant PECO-CEI from the Ministry of Higher Education
of France}

\address{Yu. G.: Algebra Section, Dept. of Mathematics,
Moscow State University, 117234 Moscow, Russia}
\email{prokhoro@mech.math.msu.su\qquad \\
prokhoro@nw.math.msu.su\\}


\title{Klein's group defines an exceptional singularity of dimension $3$}

\begin{document}
\maketitle

\section*{Introduction} 
The aim of this paper is to construct examples of 
canonical exceptional singularities. Canonical (as well as 
terminal, log terminal
and log canonical) singularities appear naturally in  the minimal model
theory and were studied by Reid, Mori, Koll\'ar, Shokurov and others. 
Recently Shokurov \cite{Sh1} introduced the notion of exceptional 
singularity, see Definition \ref{def-exc}. This notion is closely connected with
the inductive approach to the classification of singularities, 
flips, divisorial
contractions, etc. The key ingredient of the inductive approach in its
modern setting is the search of {\it complements}, that is, of good
divisors in the multiple anticanonical systems (see Definition \ref{def-compl}).
In fact the main result of \cite{Sh1} is that 
two-dimensional complements can be divided into two parts: 
regular and exceptional. Regular ones occur in $1$-, $2$-, $3$-, $4$- or 
$6$-uple anticanonical systems and have a rather simple structure. 
Exceptional ones can only occur in the neighborhood of
an exceptional singularity; they are more complicated to study but they  
belong, up to birational isomorphisms, to a
finite number of families and, at least in principle, can 
be classified \cite{Sh1}.  
By using standard arguments with
Kawamata-Viehweg vanishing and the inversion of adjunction 
(see \cite[19.6]{Ut}) these results can be applied to study 
three-dimensional log canonical singularities and, even more generally,
extremal contractions \cite[\S 7]{Sh1}.
\par
In dimension 3, 
  Shokurov \cite[\S 7]{Sh1} gave an
example of a log Del Pezzo surface with no numerical obstructions
to the existence of its blow down to an exceptional 
canonical singularity, but at the moment it is not clear
whether such a singularity really exists.
In the present paper, we construct first examples of 3-dimensional 
canonical exceptional singularities.
\par
Our examples belong to the class of quotient singularities: they are
quotients of $\CC^3$ by the action of Klein's simple group J of order
$168$, and of its central extension $J'$ by the 3-rd roots of unity, of
order 504. Remark, that these quotient singularities have been already 
investigated in relation to another problem in \cite{Mar} and \cite{Ro}.
It was proved, that they and, more generally, any quotient of $\CC^3$
by a finite subgroup of the special linear group, are crepant,
that is, admit resolutions with trivial canonical class. Thus, they are
in a sense opposite to terminal, or totally discrepant
singularities inside the class of
canonical ones. By Example \ref{primerchik}, the terminal singularities
are always nonexceptional,
so, we have a reason to start the search of
exceptional ones among crepant quotients.

We use the Miller--Blichfeldt--Dickson classification \cite{MBD}
of finite subgroups of $\SL$, and show, that for some classes
of this classification the quotient is not exceptional. These include
the reducible, imprimitive groups, the icosahedral group $H$ and its
central extension $H'$ of order 180. All these groups have a semiinvariant
of degree $\leq 3$. It is plausible, that the quotients by the subgroups
without semiinvariants of degree $\leq 3$ are exceptional. But the proof
of the exceptionality is not so easy as that of the nonexceptionality:
in the first case, one has to verify the minimal discrepancies
for all possible boundaries, whereas, in the second case, it
suffices to find one nonexceptional boundary. We prove 
the exceptionality of our examples by the rule
of contraries, in using the existence of a nonexceptional 
$1$-, $2$-, $3$-, $4$- or $6$-complement
of the canonical divisor  (according to \cite{Sh1}) 
for a nonexceptional singularity.
This implies that it suffices to verify only the boundaries
given by the semiinvariants of degree $\leq 18$. The rest of the proof
relies heavily upon the classical results of Klein on the group $J$:
enumeration of the orbits, list of invariants, the configuration of
the non-free locus of its action on the projective plane.

We will describe now briefly the contents of the article by sections.

Section 1 is preliminary, it contains definitions and some 
facts for later use. In Section 2 we describe our approach
to the prooof of the exceptionality of quotient singularities. 
In particular we prove (Corollary \ref{cor-pr}) that a quotient singularity 
$\CC^3/G$ can be exceptional only if $G$ is primitive.
In Section 3 the main result (Theorem \ref{main-th}) is proved.

\subsubsection*{Acknowledgements}
The authors would like to thank V. V. Shokurov
for useful remarks.

\section{Definitions and preliminary results}

We follow essentially the terminology and notation of \cite{Ut},
\cite{Sh} and \cite{Sh1} (see also \cite{KoP} for a nice introduction
to the subject).

\begin{definition}
Let $(X\ni P)$ be a normal singularity (not necessarily isolated) and 
let $D=\sum d_iD_i$ be a divisor on $X$ with real coefficients.
$D$ is called a boundary if $0\leq d_i\leq 1$ for all $i$.
It is called a subboundary, if it is  majorated by a boundary.
A proper birational morphism $f \colon Y\lra X$ is called a log resolution of
$(X,D)$ at $P$, if $Y$ is nonsingular near $f^{-1}(P)$ and 
$\Supp (D)\cup E$ is a normal crossing divisor
on $Y$  near $f^{-1}(P)$, where $D$ is used
to denote both the subboundary on $X$ and its proper transform on $Y$,
and $E=\cup E_i$ is the exceptional divisor of $f$.
The pair $(X,D)$ or, by abuse of language, the divisor $K_X+D$
is called {\it terminal, canonical, Kawamata log terminal (klt),
 purely log terminal
(plt),}  and, respectively,  {\it log canonical (lc)} near $P$, 
if the following
conditions are verified:

(i) $K_X+D$ is $\RR$-Cartier.

(ii) Let us write for any proper birational morphism $f \colon  Y\lra X$
$$
K_Y\equiv f^\ast (K_X+D)+\sum a(E,X,D)E ,
$$
where $E$ runs over prime divisors on $Y$, $a(E,X,D)\in\RR$,
and $a(D_i,X,D)=-d_i$ for each component $D_i$ of $D$. Then, for some
log resolution of $(X,D)$ at $P$ and for all prime
divisors $E$ on $Y$ near $P$, we have:
$a(E,X,D)>0$ (for terminal), $a(E,X,D)\geq 0$ (for canonical), 
$a(E,X,D)> -1$ and no $d_i=1$ (for klt), 
$a(E,X,D)> -1$ (for plt, without any restriction on the subboundary $D$), and,
respectively, $a(E,X,D)\geq -1$ (for lc).
\end{definition}

The coefficients $a(E,X,D)$ are called {\it discrepancies} of $f$, or of
$(X,D)$; they depend only on the discrete valuations of the
function field of $X$ associated to the prime divisors $E$, 
and not on the choice of $f$. We will
identify prime divisors with corresponding discrete
valuations, when speaking about `divisors $E$ over $X$'
without indicating, on which birational model $E$ is realized. The conditions 
given by inequalities in part (ii) of the above definition do not depend
on the choice of a log resolution. The lc (as well as terminal, canonical, 
klt, plt) condition is obviously
monotonic on $D$: if $(X,D)$ is lc at $P$, then $(X,D')$ is also lc at $P$ for
any $D'\leq D$ such that $K_X+D'$ is $\RR$-Cartier. 
Thus, for any boundary $D$, which is  a $\RR$-Cartier divisor on a
$\QQ$-Gorenstein variety we can define the {\it log canonical
threshold} of $(X,D)$ at $P$:
$$
c_P(X,D)=\max\{\alpha\in\RR \; |\; (X,\alpha D)\;\;\;
\mbox{is log canonical}\}
$$

\begin{definition}[{\cite[1.5]{Sh1}}]
\label{def-exc}
Let $(X\ni P)$ be a normal singularity  and 
let $D=\sum d_iD_i$ be a 
boundary on $X$ such that $K_X+D$ is log canonical. The pair $(X,D)$ is 
said to be {\it exceptional} if there exists at most one exceptional divisor 
$E$ over $X$ with discrepancy $a(E,X,D)=-1$. The singularity $(X,P)$ is said to
be {\it exceptional} if $(X,D)$ is exceptional for any $D$ whenewer $K_X+D$
is log canonical. 
\end{definition}

\begin{remark}\label{rem-connected}
By the connectedness result \cite[5.7]{Sh}, \cite[17.4]{Ut}, the set
of divisors with discrepancies $\le -1$ in $Y$ for any 
proper birational morphism $f \colon  Y\lra X$ is connected. Therefore,
if a log canonical pair $(X,D)$ is nonexceptional, then 
there exist infinitely many divisors with discrepancy $-1$, which can
be constructed by blowing up the compoments of intersections of pairs
of such divisors. 
\end{remark}

\begin{definition}[{\cite[5.1]{Sh}}]
\label{def-compl}
Let $(X,P)$ be a normal singularity, $D=S+B$ a subboundary, 
such that $B,S$ have no common
components, $S$ is a reduced divisor, and $B=\sum b_iB_i$ with all
$b_i< 1$, that is, $\lfloor B\rfloor =0$, where
$\lfloor\cdot\rfloor$ denotes the integer part. 
Then one sais that $K_X+D$ is {\it $n$-complemented}, if there
exists a $\QQ$-divisor $D^+$, such that the following conditions are
verified:

(i) $nD^+$ has integer coefficients ;

(ii) $nD^+\sim -nK_X$;

(iii) $nD^+-nS-\lfloor (n+1)B\rfloor\geq 0$;

(iv) $K_X+D^+$ is lc.

The divisor $K_X+D^+$ is called a {\it $n$-complement} of $K_X+D$.
Remark, that if $D$ is a boundary, then so is $D^+$.
\end{definition}

\begin{example}[{\cite[5.2.3, 5.6]{Sh}, \cite[1.5]{Sh1}}]
Let $(X,P)$ be a two-dimensional quotient singularity.
Then the following conditions are equivalent:
\begin{enumerate}
\renewcommand\labelenumi{(\roman{enumi})}
\item
$(X\ni P)$ is exceptional,
\item
$(X\ni P)$ is of type $E_6$, $E_7$ or $E_8$ (in the generalized sense
of \cite{I}, \cite{Br}); this means that the dual (weighted) graph of 
the minimal resolution of
such a singularity has a single 3-valent vertex 
$\stackrel{-b}{\circ}$, $b\ge 2$ with three 
chains issued from it, exactly one of them being
of type $\pal\stackrel{-2}{\circ}$ (see \cite{Br}, \cite{I},
 or \cite[Ch. 3]{Ut}
for a more precise description),
\item
there are no $1$- or $2$-complements $K_X+D$ such that $P\in \Supp(D)$
(but there is a $3$-, $4$- or $6$-complement),
\item
$(X\ni P)$ is analytically isomorphic to a quotient $\CC^2/G$,
where $G$ is a finite subgroup of $\GL$ without reflections of
 dihedral, tetrahedral or
icosahedral type  \cite{Br}. This means that 
the image of $G$ in $PGL_2(\CC)\simeq SO_3(\CC)$ is the dihedral, 
tetrahedral or
icosahedral group in the usual sense. 
\end{enumerate}
\end{example}

The following theorem is a consequence of the proof  of 
Shokurov's theorem \cite[7.1]{Sh1}.

\begin{theorem}
\label{Shokurov}
Let $(X\ni P)$ be a nonexceptional three-dimensional log canonical 
singularity. Then $K_X$ is either $1$-, $2$-, $3$-, $4$- or $6$-complemented.
Moreover, there exists  such a nonexceptional complement $K_X+D$.
\end{theorem} 

\begin{lemma}
\label{red}
Assume that there exists a reduced divisor $S=\sum S_i$ passing through $P$ 
such that $K_X+S$ is log canonical. Then $(X\ni P)$ is nonexceptional. 
\end{lemma} 
\begin{proof}
Take a general hyperplane section through $P$. Then, for some 
$0\le \alpha \le 1$, the log divisor $K_X+S+\alpha H$ is log canonical, but
not purely log terminal. By Remark \ref{rem-connected}, the set of divisors 
with discrepancy $a(\phantom{E},X,S+\alpha H)=-1$ on a resolution $Y\to X$
is connected, and by construction, we have 
at least two of them, one coming from $S$ and the other one, say $E$,
exceptional over $X$. So
we can get infinitely many divisors with discrepancy $-1$ in blowing up 
the curves of intersections.
\end{proof}

\begin{example}
By definition and Lemma \ref{red} any  three-dimensional cDV-singularity
is nonexceptional.
\end{example}

\begin{example}
\label{primerchik}
Let $(X\ni P)$ be a three-dimensional terminal singularity and
$S\in |-K_X|$ a general element 
(in the Gorenstein case, we should additionally suppose that $S\ni P$). 
By \cite[6.4]{YPG} and the inversion of adjunction 
\cite[3.3, 3.12, 5.13]{Sh}, \cite[17.6]{Ut}
$K_X+S$ is purely log terminal (and even canonical).
Therefore, all terminal singularities are nonexceptional.
Of course, these arguments use the classification of 
terminal singularities. Shokurov  (cf. \cite[6.5]{YPG})
posed the problem to prove this fact directly.
\end{example}

\section{Quotient singularities}

Now let $(X\ni P)$ be a three-dimensional quotient singularity, i.~e.
$(X\ni P)=(\CC^3\ni 0)/G$, where $G\subset \GL$ is a finite subgroup.
We may assume that $G$ contains no quasi-reflections, for if 
$G$  contains quasi-reflections, then there exists another subgroup
$G'$ of $\GL$, which does not contain a quasi-reflection and such that
$\CC [y_1,y_2,y_3]^G\simeq\CC [y_1,y_2,y_3]^{G'}$ as $\CC$-algebras.
\par
Let $\pi\colon V\to X$ be the quotient morphism, where $V=\CC^3$.
Let $D$ be a boundary on $X$ and let $D':=\pi^*D$. By \cite[2.2]{Sh},
\cite[20.3]{Ut}, \cite[]{KoP} $K_X+D$ is log canonical 
(resp., plt, klt) iff so 
is $K_{V}+D'$. 

\begin{lemma}\label{quot-lem1}
In the above notation,
$(X,D)$ is exceptional iff $(V,D')$ is.  
\end{lemma}
\begin{proof}
Assume that $(V,D')$ is exceptional. Then by \cite[3.1]{Sh2},
\cite[17.10]{Ut} there exists a blow-up $g\colon W\to V$
such that $K_{W}+S+B'={g}^{*}(K_{V}+D')$ is purely log terminal, where
$S$ is the (irreducible) exceptional divisor of $g$ and $B'$ is the proper 
transform of $D'$. This blow-up is unique up to isomorphism, because
$g$ is projective and $\rho(W/V)=1$ \cite[6.2]{Ut}.  So we can
define an action of $G$ on $W$ making $g$ equivariant. 
Let $\varphi\colon W\to Y$
be the  quotient morphism and put
 $E:=\varphi(S)$ and $B:=\varphi(B')$. 
We have the following commutative diagram
\begin{equation}
\label{1}
\begin{CD}
W@>{\varphi}>>Y\\
@V{g}VV@VfVV\\
V@>{\pi}>>X.\\
\end{CD}
\end{equation}
By the ramification formula
\begin{equation}
\label{2}
K_{W}+S+B'=\varphi^*(K_{Y}+E+B),
\end{equation}
whence
\begin{equation}
\label{3}
K_{Y}+E+B=f^*(K_X+D).
\end{equation}
By \cite[2.2]{Sh},
 \cite[20.3]{Ut} $K_{Y}+E+B$ is purely log terminal.
Therefore, $K_X+D$ is exceptional. 
\par
Conversely, assume that $K_X+D$ is exceptional. Let $W$
be the normalization of $V\times_{X}Y$. We again have the 
diagram (\ref{1}) and relations (\ref{2}), (\ref{3}).
Similarly, by \cite[2.2]{Sh}, \cite[20.3]{Ut},
$K_{W}+S+B'$ is purely 
log terminal, hence  $K_{V}+D'$ is exceptional. 
\end{proof}

\begin{lemma}\label{quot-lem2}
In the notation of Lemma \ref{quot-lem1}, assume that $G$ 
has a semiinvariant of degree
$\le 3$. Then $(X\ni P)$ is nonexceptional.
\end{lemma}

\begin{proof}
Let $\psi$ be  such a semiinvariant of minimal degree $d\le 3$, 
let $D'$ be its zero locus, and let $D:=\pi(D')_{\red}$. Then $D'=\pi^*(D)$.
By Lemma \ref{red} it is sufficient to show that  $K_{V}+D'$
is log canonical. Note that $\psi$ is homogeneous, so $D'$ is a cone over
a plane curve $C$ of degree $d\le 3$. 
 If  $C$ is nonsingular, 
then $D'$ has a unique singular point which can be resolved by
only one blow-up with $a(\nu, V,D')=2-\mult_0\psi\ge -1$,
so  $K_{V}+D'$
is log canonical in this case.
  If $C$ has a unique singular point, then 
$G$ has an eigenvector, say $v\in \CC^3$. But then $G$ has an 
invariant plane, orthogonal to $v$ and we can take $d=1$ and $C$
is nonsingular. 
The same arguments work if
$C$ has exactly two singular points. In the remaining case  $C$ 
is the union of three lines in general 
position. Then $(V,D')$ is already log nonsingular
(that is, may be taken as its own log resolution), and so $K_{V}+D'$
is log canonical.      
\end{proof}

The subgroup $G\in\GL$ is said to be {\it reducible}, if it has a proper
invariant subspace in $\CC^3$. $G$ is called {\it imprimitive}, if there exists
a triple of lines $L_1,L_2,L_3$ in $\CC^3$, 
permuted by $G$ via a representation
$G\lra S_3$, where $S_n$ denotes the symmetric group of permutations
on $n$ elements.

\begin{corollary}
\label{cor-pr}
If  $(X\ni P)$ is  exceptional, then $G$ is irreducible 
and primitive.
\end{corollary}

\begin{proof}
If $G$ is reducible, it has a semiinvariant of degree 1. If it
is imprimitive, it has a semiinvariant of degree 3, defining the
union of three planes spanned by pairs of lines from $L_1,L_2,L_3$.
The result follows by Lemma \ref{quot-lem2}.
\end{proof}

The finite subgroups of $G\subset\GL$ were classified by
Miller--Blichfeldt--Dickson \cite{MBD} modulo extension by
certain scalar matrices (compare with \cite{P}). There are 9 types
of such groups, denoted by A,B, \ldots ,J in \cite{MBD}.
As soon as we are looking for those which yield exceptional 
quotient singularities, we have to test only primitive
irreducible ones. They belong to the $6$ types E, F, G, H, I, J.
The orders of the associated collineation groups 
$PG=G/(G\cap\CC^\ast)\subset\PGL$ 
are $36$, $72$, $216$, $60$, $360$, $168$; 
the first three are solvable,
and the last three are simple. The collineation 
groups from G to J have their names:
the Hessian group, the icosahedral one, the alternating group 
of degree $6$, and, finally,
Klein's simple group. 

\begin{prop}
The quotients of $\CC^3$ by the subgroups of $\SL$
 of type H are nonexceptional.
\end{prop}
\begin{proof}
There are two such groups 
(see, e. g., \cite{P} or \cite{Ro}): the icosahedral group $G$
of order $60$, and its central extension $G'$ of order $180$.
The invariants of $G$ are well-known; see, for example,
\cite[Sect. 116]{MBD}. 
There is an invariant of degree $2$, 
which is a semiinvariant of $G'$.
(This also follows from the fact that $G$ is a subgroup of
$SL_3(\RR)$ and therefore has an invariant quadratic form 
\cite[4.2.15]{Sp}).
The result follows by Lemma \ref{quot-lem2}.
\end{proof}
\subsection{}
\label{logic}
Now we will explain the logic of our approach to the proof
of the exceptionality of a quotient singularity, which
will be applied in the next section to Klein's group.
Assume that $(X\ni P)$ is nonexceptional. By Theorem
\ref{Shokurov} there exists a nonexceptional log canonical
$K_X+D$ such that $n(K_X+D)\sim 0$ for $n\in\{1, 2, 3, 4, 6\}$.
Further, we will use notations of Lemma \ref{quot-lem1}. The integer divisor
$F:=nD'$ locally near $0$ can be defined by a 
semiinvariant function, say $\psi$.
\par
Moreover $nD\sim -nK_X$ iff the form  $\psi(dx_1\wedge dx_2\wedge dx_3)^{-n}$
is invariant, i.~e. 
\begin{equation}\label{character}
g(\psi)=\det(g)^n\psi\qquad\text{ for all}\qquad g\in G.
\end{equation}
\par  
Denote $d:=\mt{mult}_0(\psi)$. 
\label{begin}
Let $\sigma\colon W\to V=\CC^3$ be the blow-up of the origin
and let $S\simeq\PP^2$ be the exceptional divisor.
Then $K_W=\sigma^*K_V+2S$ and
$\sigma^*F=R+dS$, where $R$ is the proper transform of $F$.
Therefore the  discrepancy of $S$ is 
$$
a(S,V,D')=2-\frac{1}{n}\mt{mult}_0(\psi)\ge -1.
$$
So we have $d=\mt{mult}_0(\psi)\le 3n\leq 18$.

Further
\begin{equation}
K_W+S+\frac{3}{d}R=\sigma^*(K_V+\frac{3}{d}F).
\end{equation}
By \cite[Lemma 3.10]{KoP} $K_V+\frac{3}{d}F$ is log canonical iff so is 
$K_W+S+\frac{3}{d}R$. In this case the pair $(V,\alpha F)$ is 
exceptional for all $0\leq\alpha\leq \frac{3}{d}$ iff 
$K_W+S+\frac{3}{d}R$ is purely log terminal; by the above, we need
this assertion only for $\alpha =\frac{1}{n}$ with $n\in\{1, 2, 3, 4, 6\}$,
but we will have a stronger property with $\alpha$ not necessarily of this
form.
The plt condition for $K_W+S+\frac{3}{d}R$ is equivalent to 
that $K_S+\frac{3}{d}C$
is Kawamata log terminal, where $C=R\cap S$ 
\cite[5.13]{Sh}, \cite[17.6]{Ut}.
It is clear that $C$ is given by the equation 
$\psi_{\min}=0$, where  $\psi_{\min}$ is 
the homogeneous component of  $\psi$ of minimal  degree $d$.
Therefore we have

\begin{proposition}
\label{check}
In the above notations, if $K_S+\frac{3}{d}C$
is Kawamata log terminal, then $(V,\alpha F)$ is exceptional for any
$0\leq\alpha \leq\frac{3}{d}$.
\end{proposition}

\begin{remark}
If $K_S+\frac{3}{d}C$ is log canonical, but not Kawamata log terminal, 
then $(V,F)$ is nonexceptional. If $K_S+\frac{3}{d}C$
is not log canonical, then we can conclude nothing (possibly in this 
situation we have to consider some weighted blow-up). 
\end{remark}

\begin{lemma}
\label{bound}
Notations as above.
Assume that $K_S+\frac{3}{d}C$
is not Kawamata log terminal and $C$ is a singular irreducible curve. 
Then there exists an orbit of $G$ consisting of at most $10$ 
singular points of $C$. 
\end{lemma}

\begin{proof}
Let $P\in C$ be a singular point of maximal multiplicity $m$
and let $r$ be the number of points in the orbit $G\cdot P$.
Denote by  $c=c(S,C)$ the (global) log canonical threshold of $(S,C)$.
By our assumption $c\le 3/d$. On the other hand, $c\ge 1/m$  
\cite[Lemma 8.10]{KoP}. This gives us $d\le 3m$. Let $g$ be the genus of the 
normalization of $C$. Then 
$$
0\le g\le \frac{(d-1)(d-2)}{2}-r\frac{m(m-1)}{2}.
$$
Taking into account that $d\le 3m$, we obtain 
$$
-2\le 2g-2\le (9-r)m(m-1).
$$
This implies the assertion.
\end{proof}

\section{Klein's group}

The aim of this section is to prove the exceptionality of the
quotients of $\CC^3$ by the subgroups of $\SL$ of type J
in the classification of \cite{MBD}. There are two such groups,
see \cite{P} or \cite{Ro}:  {\it Klein's simple group} $J_{168}$ 
of order $168$ and its central 
extension $J'_{504}$ of order $504$.
\begin{theorem}\label{main-th}
Let $G\subset\SL$ be $J_{168}$ or $J'_{504}$. Then the singularity 
of the quotient $\CC^3/G$ 
at the origin is exceptional. 
\end{theorem}

We will briefly describe the irreducible representation of $J_{168}$ in
$\CC^3$ following \cite{W} and \cite{Kl}. 
Another description of $J_{168}$ see in \cite[Sect. 4.6]{Sp}. 
Let $y_1,y_2,y_3$ be coordinates in
$\CC^3$. The group $J_{168}$ is generated by $3$ elements $\tau ,\chi ,
\omega$ of orders $7$, $3$, $2$ respectively, and the representation is
defined by
\[ \begin{array}{lcl}
\tau =\left(\begin{array}{ccc}
\epsilon &  & 0 \\
 & \epsilon^{2} &  \\
0 &  & \epsilon^{4} \end{array}\right) & , &
\epsilon =\exp\left(\frac{2\pi i}{7}\right) ,\vspace*{2 ex} \\
\chi =\left(\begin{array}{ccc}
0 & 0 & 1 \\
1 & 0 & 0 \\
0 & 1 & 0 \end{array}\right) & , & 
\omega =\left(\begin{array}{ccc}
\alpha & \beta & \gamma \\
\beta & \gamma & \alpha \\
\gamma & \alpha & \beta \end{array}\right) , \vspace*{2.5 ex}\\
\alpha =-\frac{2\sin\frac{8\pi}{7}}{\sqrt{7}}\; ,\;\; &
\beta =-\frac{2\sin\frac{4\pi}{7}}{\sqrt{7}}\; ,\;\; &
\gamma =-\frac{2\sin\frac{2\pi}{7}}{\sqrt{7}},
\end{array}\]

The second group $J'_{504}$ is generated by $J_{168}$ and the
scalar matrix with $\exp \frac{2\pi i}{3}$ on the diagonal.

Now, we will describe the semiinvariants of these groups. First of
all, since $J_{168}$ is simple, all its semiinvariants are indeed invariants.
They are also semiinvariants of $J'_{504}$. According to \cite{Kl}
(see also \cite{W}), the algebra of invariants $A=\CC [y_1,y_2,y_3]^{J_{168}}$
is generated by four homogeneous polynomials $f,\Delta ,\CCC,\KKK$ of degrees
$4$, $6$, $14$, $21$ respectively, with one basic relation between them:

\begin{eqnarray*}
f & = & y_{1}^{3}y_{3}+y_{2}^{3}y_{1}+y_{3}^{3}y_{2},\vspace*{2 ex} \\
\Delta & = & \frac{1}{54}\mbox{Hess}\, (f), \vspace*{2 ex}\\
\CCC & = & \frac{1}{9}\left|\begin{array}{cccc}
f_{y_{1}y_{1}}^{\prime\prime} &
f_{y_{1}y_{2}}^{\prime\prime} &
f_{y_{1}y_{3}}^{\prime\prime} &
\Delta_{y_{1}}^{\prime} \\
f_{y_{2}y_{1}}^{\prime\prime} &
f_{y_{2}y_{2}}^{\prime\prime} &
f_{y_{2}y_{3}}^{\prime\prime} &
\Delta_{y_{2}}^{\prime} \\
f_{y_{3}y_{1}}^{\prime\prime} &
f_{y_{3}y_{2}}^{\prime\prime} &
f_{y_{3}y_{3}}^{\prime\prime} &
\Delta_{y_{3}}^{\prime} \\
\Delta_{y_{1}}^{\prime} & \Delta_{y_{2}}^{\prime} & 
\Delta_{y_{3}}^{\prime} & 0
\end{array}\right| ,\end{eqnarray*}
\begin{eqnarray*}
\\
\KKK & = & \frac{1}{14}\left|\begin{array}{ccc}
f_{y_{1}}^{\prime} & \Delta_{y_{1}}^{\prime} & \CCC_{y_{1}}^{\prime} \\
f_{y_{2}}^{\prime} & \Delta_{y_{2}}^{\prime} & \CCC_{y_{2}}^{\prime} \\
f_{y_{3}}^{\prime} & \Delta_{y_{3}}^{\prime} & \CCC_{y_{3}}^{\prime} 
\end{array}\right| . \vspace*{4 ex}\end{eqnarray*}

\begin{eqnarray}
\KKK^{2} & = & \CCC^{3}+1728\Delta^{7}+1008\CCC\Delta^{4}f-88\CCC^{2}\Delta f^{2}
-60032\Delta^{5}f^{3} \nonumber\\
 & & \rule{2 em}{0 pt} +1088\CCC\Delta^{2}f^{4}+22016\Delta^{3}f^{6}
-256\CCC f^{7}-2048\Delta f^{9}.
\label{Klein-1}\end{eqnarray}

In the proof of the Theorem, we will follow the logic scheme 
outlined in Subsect. \ref{logic}. In the notations of \ref{logic} let $F$
be given by the equation $\psi =0$, where $\psi = \sum a_{ijkl}f^i
\Delta^j\CCC^k\KKK^l\;\; (i\geq 0,j\geq 0,k\geq 0,l=0$ or 1, such that the
the lowest degree of non-zero terms $d=\mbox{mult}_0\psi = \deg\psi_{\min}
\leq 18$. For $G=J_{168}$, up to a constant factor of 
proportionality, $\psi_{\min}$
can be only one of the following functions:
$$
f,\De , f^2, f\De , \la f^3+\mu\De^2, \la f^2\De+\mu \CCC,
\la f^4+\mu \De^2f, \la\De^3+\mu f^3\De +\nu f\CCC,
$$
where $\la ,\mu ,\nu$ are arbitrary complex constants. 
The list of initial forms for $G=J'_{504}$ is even shorter: 
according to (\ref{character}), the character of $\psi$
is equal to $g\mapsto(\det g)^n$, hence it is trivial,
because $G\subset\SL$. So, only polynomials of degree
divisible by 3 should be kept in the case of  $J'_{504}$.
The rest of the proof depends only on the pair $(S,C)=(\PP^2,(\psi_{\min}))$, 
and is done simultaneously for the two groups.\smallskip

\subsection*{Case 1} $\psi_{\min}=f^k$ or $\De^k$; $\deg\psi_{\min}= 4k$ or, 
resp.
$6k$. The reduced curve $C_{{\mbox{{\tiny red}}}}$ is nonsingular, so
the pair $(S,C_{\mbox{{\tiny red}}})$ is log nonsingular. 
Hence  $(S,\alp C_{\mbox{{\tiny red}}})$ is klt for
any $\alp <1$. Hence $(S,\frac{3}{d}C)=(S,\frac{3k}{d}C_{\mbox{{\tiny red}}})$
is klt. By Proposition \ref{check}, the pair $(V,\alpha F)$ 
is exceptional for any $0\leq\alpha \leq\frac{3}{d}$, and we are done.
\smallskip

\subsection*{Case 2} $\psi_{\min}=f^i\De^j, \; 10\leq d=4i+6j\leq 18$. Here
$C_{\mbox{{\tiny red}}}=\{ f\De =0\}$ is singular only at the points
of intersection of Klein's quartic $C_1=\{ f=0\}$ and 
its Hessian curve $C_2=\{ \De =0\}$,
that is at the inflection points of Klein's curve. It is known that they
form one orbit of length 24 under the action of $J_{168}$ with representative
$(1:0:0)$. Hence they are ordinary double points of $C_{\mbox{{\tiny red}}}$,
and thus, $(S,C_{\mbox{{\tiny red}}})$ is log nonsingular. Hence
$(S,\alp C_1+\beta C_2)$ is klt for any $\alp <1,\beta < 1$. Hence
$(S,\frac{3}{d}(iC_1+jC_2))$ is klt, and we are done.\smallskip

Cases 1 and 2 cover all the invariant curves of degree $<12$, as well as
multiples of $C_1$ and $C_2$.\smallskip

\subsection*{Case 3} $C$ is a reduced irreducible curve of degree $d\geq 12$.
If it is nonsingular, we are done, because $(S,C)$ is log nonsingular,
and $\frac{3}{d}<1$. Assume that $(S,\frac{3}{d}C)$ is not klt. Then $C$ 
cannot be singular by Lemma \ref{bound}, because
$J_{168}$ does not have orbits in $\PP^2$ of length $\leq 10$. For the orbits of
Klein's group, see, e.~g. \cite[Sect. 120]{W}; the possible
lengths are $21$, $24$, $28$, $42$, $56$, $84$, $168$. 
Hence $C$ is nonsingular,
and this is a contradiction. Hence $(S,\frac{3}{d}C)$ is always klt
in this case, and we are done.\smallskip

\subsection*{Case 4} $C$ is reducible or non-reduced of degree $12\leq d\leq 
18$.
The irreducible components of $C$ are permuted by the action of $J_{168}$.
The length of the orbit (if not 1, for an invariant component), is
$\geq 7$, for $J_{168}$ has no non-trivial homomorphisms to symmetric
groups $S_p$ with $p<7$. Moreover, there are no straight lines as irreducible
components, because the action on the dual projective plane of lines in
$\PP^2$ is the same as on $\PP^2$ itself, and hence the minimal orbit
length is 21 (there is indeed an invariant of degree 21 which factors
into the product of linear forms, namely, $\KKK$). So, discarding the multiples
or combinations of $C_1$ and $C_2$ covered by Cases 1 and 2,
we have three subcases: A. $f\not\mid
\psi_{\min}$ and $\De\not |
\psi_{\min}$. Then $d=$ $14$, $16$, or $18$, and $C$ is the union of $d/2$
conics. B. $f|\psi_{\min}$. Then $d=16$ or $18$, and $\psi_{\min}=f
(\la f^3+\mu\De^2)$ with $\la\mu\neq 0$ or, resp., $\psi_{\min}=f
(\la f^2\De +\mu C)$ with $\mu\neq 0$. C. $\De |\psi_{\min}$. Then
 $\psi_{\min}=\De (\la f^3+\mu\De^2)$ with $\la\mu\neq 0$.
\smallskip

\subsubsection*{Subcase A} $d=18$ can be eliminated, because 9$\not |$168.
If $d=16$, then the stabilizer of any conic component $\G$ of $C$ is of
order 21. Hence, the orbit of the generic point of $\G$ is of length
$21$, which is impossible, because there is only one orbit of length
$21$ in $\PP^2$: `die achtz\"{a}hlige Pole' in the classical terminology.
A similar argument shows that $d=14$ is also impossible.
\smallskip

\subsubsection*{Subcase B} If $d=16$, then $C=C_1\cup C'$, where $C_1$ is
Klein's curve, and $C'=\{ \la f^3+\mu\De^2=0\}$ is an irreducible
curve of degree $12$. The $24$ points $Q_i$ of $C_1\cap C'$ are ordinary cusps
of $C'$, and the local indices $(C_1\cdot C')_{Q_i}=2$. By Bezout
Theorem, these are the only intersection points. Assume that there is a
singular point $R\in C'$, different from $Q_i$. Let $\de$ be the length
of the orbit of $R$, and $m'$ its multiplicity. Then $p_a(C')=55\geq
24+\de\frac{m'(m'-1)}{2}$. Taking into the account the possible lengths
of orbits, and the uniqueness of orbits of any length $<84$, we obtain
two possible values $\de = 21$ or 28, and $m'=2$. (Remark, that
$\de =28$ really occurs for the dual of Klein's curve, and the 
corresponding singular points are ordinary double). 
So, the maximal multiplicity
of singular points of $C$ is  $m=3$, attained at the points $Q_i$.
By \cite[Lemma 8.10]{KoP}, we have the following estimate for the
log canonical threshhold: $c_Q(S,C)\geq\frac{1}{\mt{mult}_QC}\geq\frac{1}{3}$.
As $\frac{3}{d}< \frac{1}{3}$ in our case, we are done.\smallskip

\noindent If $d=18$, then $C=C_1\cup C''$, where $C''=\{ \la f^2\De +\mu C=0\}$
is irreducible of degree 14. Similarly to the 
above, we see that the only possibility
for the intersection locus $C_1\cap C''$ is the orbit of the point $\{ y_1=
\exp \frac{2\pi i}{3}y_2=\exp \frac{4\pi i}{3}y_3\}$ of length 56 (the 
other candidate, the orbit of lenght 28 taken with multiplicity 2, is
eliminated because it is not contained in $C_1$!). So, the intersections are
transversal. Assuming that there is an extra singular point of $C$,
we find at once that the length of its orbit should be 21, and in this case
the genus of $C''$ is 1. Klein's group cannot act non-trivially
on an elliptic curve. So, $C$ has only 56 ordinary double points
as singularities, and hence $(S,\alp C)$ is klt for all $\alp <1$.
\smallskip

\subsubsection*{Subcase C}
 $C=C_2\cup C'$, where $C_2=\{\De =0\}$ is
the Hessian curve, and $C'=\{ \la f^3+\mu\De^2=0\}$ is the irreducible
curve of degree 12 from Subcase B. The 24 points $Q_i$ are the points
of triple intersection of $C_2,C'$, so by Bezout, there are no other
points of intersection. By the argument of Subcase B, $C'$ cannot
acquire singularities, worse than double points. So, the points $Q_i$
are the singular points of $C$ of maximal multiplicity $m=3$.
The same argument as in Subcase B ends the proof.\smallskip

The cases 1-3 cover all possible invariant curves of degree $\leq 18$,
and in all these cases the pair $(V,\alp F)$ is exceptional as soon as 
it is log canonical. By \ref{logic}, this ends the proof of the Theorem.

\begin{remark}
The exceptional divisor $E$ of $f\colon Y\to X=\CC^3/G$ of Sect. 2
(see (\ref{1}))
for Klein's group $G$ is the weighted projective plane 
$\PP(4,6,14)=\PP^2/G$. The different $\Diff_E(0)$ (see \cite[Ch. 1]{Sh1}
for the definition) in this case is 
$(1/2)\Gamma$, where $\Gamma$ is an irreducible curve, the image of 
$21$ lines of fixed points of the elements of order $2$ in $G$.
The curve $\Gamma$ is given by the equation
\begin{eqnarray*}
0 & = & \CCC^{3}+1728\Delta^{7}+1008\CCC\Delta^{4}f-88\CCC^{2}\Delta f^{2}
-60032\Delta^{5}f^{3} \nonumber\\
 & & \rule{2 em}{0 pt} +1088\CCC\Delta^{2}f^{4}+22016\Delta^{3}f^{6}
-256\CCC f^{7}-2048\Delta f^{9}
\end{eqnarray*}
of weighted degree $42$ (see (\ref{Klein-1})). 
It is clear that $(E, \Diff_E(0))$ is a log Del Pezzo surface.
The weighted projective plane $E=\PP(4,6,14)$ has three singular points 
of types $A_1$, $A_2$ and $\CC^2/\ZZ_7(2,3)$.
The curve $\Gamma$ passes through the first of them and
has two singular points: a simple cusp and a tacnode point
(lying in the non singular part of $E$).
It is easy to compute that $K_E+(1/2)\Gamma$ is 
$1/7$-log terminal, so Shokurov's invariant $\delta$ (see \cite{Sh1})
is $0$ in our case.
A  very interesting question is to compute the minimal
complement of $K_E+\Diff_E(0)$.
\end{remark}

\end{document}